\documentclass[12pt]{amsart}

\makeatletter

\usepackage{amsmath,amssymb,doublespace,euscript}

\newcommand{\LyX}{L\kern-.1667em\lower.25em\hbox{Y}\kern-.125emX\spacefactor1000}

\theoremstyle{plain}    
\newtheorem{thm}{Theorem}
\numberwithin{figure}{section} 
\theoremstyle{plain}    
\newtheorem{lemma}[thm]{Lemma} 

\theoremstyle{remark}    
\newtheorem{case}{Case}[thm]

\makeatother

\newcount\theTime
\newcount\theHour
\newcount\theMinute
\newcount\theMinuteTens
\newcount\theScratch
\theTime=\number\time
\theHour=\theTime
\divide\theHour by 60
\theScratch=\theHour
\multiply\theScratch by 60
\theMinute=\theTime
\advance\theMinute by -\theScratch
\theMinuteTens=\theMinute
\divide\theMinuteTens by 10
\theScratch=\theMinuteTens
\multiply\theScratch by 10
\advance\theMinute by -\theScratch

\def\today{{\number\day\space
 \ifcase\month\or
  January\or February\or March\or April\or May\or June\or
  July\or August\or September\or October\or November\or December\fi
 \space\number\year}}

\newcommand\Ahat{{\widehat A}}

\newcommand\Bc{{\mathcal{B}}}

\newcommand\Cpx{{\mathbf C}}

\newcommand\freeprod{\operatornamewithlimits{\ast}}

\newcommand\HEu{{\EuScript H}}                   

\newcommand\id{\mathrm{id}}

\newcommand\Ints{{\mathbf Z}}

\newcommand\KEu{{\EuScript K}}                   

\newcommand\LEu{{\EuScript L}}                   

\newcommand\Nats{{\mathbf N}}

\newcommand\otdt{\otimes\cdots\otimes}

\newcommand\oup{^{\mathrm o}}

\newcommand\REu{{\EuScript R}}                   

\begin{document}

\begin{spacing}{1.2}

\title{Free products of exact groups}

\date{2 September, 1999}
 
\author{Kenneth J.\ Dykema}

\begin{abstract}
It has recently been proved that the class of unital exact C$^*$--algebras is closed
under taking reduced amalgamated free products.
Here the proof is presented of a special case: that the class of exact discrete groups
is closed under taking free products (with amalgamation over the identity element).
The proof of this special case is considerably simpler than in full generality.
\end{abstract}

\maketitle

\markboth{}{}
 
\section*{Introduction.}

A C$^*$--algebra $A$ is said to be {\em exact} if taking the minimal tensor product with $A$
preserves exactness of short exact sequences;
namely, $A$ is exact if and only if whenever
\begin{equation*}
0\to I\hookrightarrow B\overset\pi\to B/I\to0
\end{equation*}
is a short exact sequence of C$^*$--algebras and $*$--homomorphisms, it follows that the resulting sequence
\begin{equation}
\label{eq-minAexact}
0\to I\otimes_{\min}A\hookrightarrow B\otimes_{\min}A\overset{\pi\otimes\id_A}\longrightarrow(B/I)\otimes_{\min}A\to0
\end{equation}
is exact.
The issue of exactness was first raised by S.~Wassermann's result~\cite{Wa76} that the full group
C$^*$--algebra of a free group is not exact.
It is known that every nuclear C$^*$--algebra $A$ is exact and that C$^*$--subalgebras of exact C$^*$--algebras are exact;
see Wassermann's monograph~\cite{Wa:SNU} for a good introduction to exact C$^*$--algebras.

Let us say that a group $G$ is exact if it's reduced group C$^*$--algebra
$C^*_r(G)$ is an exact C$^*$--algebra, where we take $G$ with the discrete topology;
by definition $C^*_r(G)\subseteq\Bc(l^2(G))$
is the C$^*$--algebra generated by the left regular representation of $G$.
(We will always take groups to have the discrete topology;
the more general case of locally compact groups has been explored in~\cite{KW:exgp}.)
Many groups are known to be exact.
For example, amenable groups are exact, because their reduced group C$^*$--algebras are nuclear.
Choi proved~\cite{Choi:Z2Z3} that $C^*_r(\Ints_2*\Ints_3)$ embeds in a nuclear C$^*$--algebra,
where $\Ints_2*\Ints_3$ is the free product of the two element group and the three element group;
this implies that $\Ints_2*\Ints_3$ is an exact group, and passing to subgroups shows that
every nonabelian free group is exact.
Hyperbolic groups are known to be exact by Adams' result~\cite{Adams:amen} implying that
their reduced group C$^*$--algebras embed in nuclear C$^*$--algebras; (see also the recent proof of Germain~\cite{Ger:exgp}).
The class of exact groups is known to be closed under passing to subgroups,
under taking direct limits and under taking direct sums;
moreover, by~\cite[Theorem 5.1]{KW:perm}, extensions of exact groups by exact groups are exact.
It is in fact conjectured that all discrete groups are exact.

In~\cite{D:ExactFP} we proved that every C$^*$--algebra constructed
as a reduced amalgamated free product of exact C$^*$--algebras is exact;
an immediate consequence is that the class of exact groups is closed under taking amalgamated free products.
Here we will give a simplified version of the argument by way of proving a special case of the result,
namely that the class of exact groups is closed under taking free products, (i.e.\ amalgamating over only the identity
elements).

\section*{Exactness of free products.}

Throughout this section, let $G=\freeprod_{\iota\in I}G_\iota$ be the free product of
a family $(G_\iota)_{\iota\in I}$ of groups.
Recall that $G$ is constructed as
\begin{equation}
\label{eq-FreeProdDef}
G=\{e\}\cup\bigl\{g_1g_2\cdots g_n
\mid g_j\in G_{\iota_j}\backslash\{e\},\,\iota_1\ne\iota_2,\iota_2\ne\iota_3,\ldots,\iota_{n-1}\ne\iota_n\bigr\},
\end{equation}
where multiplication is performed by concatenation and reducing.
Our goal is to show that $G$ is exact under the hypothesis that all the $G_\iota$ are exact.
The proof has two main steps, namely Lemmas~\ref{lem-stepI} and~\ref{lem-stepII} below; note that in Lemma~\ref{lem-stepI}
the $G_\iota$ are not required to be exact.
We shall also use the following result.
\begin{lemma}
\label{lem-W} Let $A$ be a C$^*$--algebra and let $\pi:A\to\Bc(\HEu)$ be a faithful $*$--representation
of $A$ on a Hilbert space $\HEu$.
Suppose that for every finite subset $\omega\subseteq A$ and every $\epsilon>0$ there is an exact C$^*$--algebra $D$
and there are completely positive contractions $\phi:A\to D$ and $\psi:D\to\Bc(\HEu)$ such that
$\|\psi\circ\phi(x)-\pi(x)\|<\epsilon$ for every $x\in\omega$. 
Then $A$ is an exact C$^*$--algebra.
\end{lemma}
\noindent
This lemma was proved by Wassermann~\cite{Wa:nuclemb}, but where the C$^*$--algebra $D$ was required to be finite dimensional.
However, the proof of the more general result stated above is only a minor modification of Wassermann's proof;
it can be found in~\cite{D:ExactFP}.

A natural length function on the free product of groups is the {\em block length}, which we will denote $\ell$;
an element $g=g_1g_2\cdots g_n$ as in~\eqref{eq-FreeProdDef} has block length $\ell(g)=n$, while $\ell(e)=0$.
Let $W_n=\{g\in G\mid\ell(g)\le n\}$, regard $l^2(W_n)$ as a subspace of $l^2(G)$ and let $P_n:l^2(G)\to l^2(W_n)$
denote the projection.
Let $\Phi_n:\Bc(l^2(G))\to\Bc(l^2(W_n))$ be the unital completely positive map $\Phi_n(x)=P_nxP_n$.
We now find maps going the other way which give in the limit an approximate inverse for $\Phi_n$ on elements of $C^*_r(G)$.
\begin{lemma}
\label{lem-stepI}
There are completely positive unital maps $\Psi_n:\Bc(l^2(W_n))\to\Bc(l^2(G))$, for all $n\in\Nats$, so that
\begin{equation}
\label{eq-PsiPhi}
\forall a\in C^*_r(G)\qquad\lim_{n\to\infty}\|\Psi_n\circ\Phi_n(a)-a\|=0.
\end{equation}
\end{lemma}
Before beginning the proof, note that it is entirely plausible that such maps $\Psi_n$ exist.
If $a\in C^*_r(G)$ and if $a$ happens to be a linear combination of left translators $\lambda_g$, where $\ell(g)\le k$,
(the set of such elements $a$, for $k$ arbitrary, is dense in $C^*_r(G)$),
then $a$ can be recovered from the compression $\Phi_n(a)$ as soon as we have $n\ge k$;
indeed, $a\delta_e\in l^2(W_n)$ and $a$ can be found from $a\delta_e$.
However, we must (approximately) recover $a$ from $\Phi_n(a)$ by use of a completely positive map.

\begin{proof}
We will define an isometry $V_n:l^2(G)\to l^2(W_n)\otimes l^2(G)$ and let $\Psi_n:\Bc(l^2(W_n))\to\Bc(l^2(G))$
be given by $\Psi_n(x)=V_n^*(x\otimes1)V_n$.
If $g\in G$ with $\ell(g)>0$, let $\iota_1,\ldots,\iota_{\ell(g)}\in I$ and $g_j\in G_{\iota_j}\backslash\{e\}$ be such that
$\iota_j\ne\iota_{j+1}$ and $g=g_1g_2\ldots g_{\ell(g)}$.
Let $m=\bigl[\frac{n+1}2\bigr]$.
Let
\begin{equation}
\label{eq-Vndef}
V_n\delta_g=
\begin{cases}
\delta_g\otimes\delta_e&\text{if }\ell(g)\le m \\[1ex]
\frac{\sqrt{n-\ell(g)}}{\sqrt{n-m}}\delta_g\otimes\delta_e
+\sum_{j=m}^{\ell(g)-1}\frac1{\sqrt{n-m}}\delta_{g_1\ldots g_j}\otimes\delta_{g_{j+1}\ldots g_{\ell(g)}}&\text{if }m\le\ell(g)\le n \\[1ex]
\sum_{j=m}^{n-1}\frac1{\sqrt{n-m}}\delta_{g_1\ldots g_j}\otimes\delta_{g_{j+1}\ldots g_{\ell(g)}}&\text{if }n\le\ell(g).
\end{cases}
\end{equation}
In order to show that~\eqref{eq-PsiPhi} holds, we may without loss of generality take $a=\lambda_h$ for $h\in G$,
because the linear span of these is dense in $C^*_r(G)$.
For every $g\in G$ we have that $\ell(g)-\ell(h)\le\ell(hg)\le\ell(g)+\ell(h)$, so taking $n$ to be very much larger than $\ell(h)$,
the following cases are the only possibilities, and we find that $\Psi_n\circ\Phi_n(\lambda_h)\delta_g=c(n,h,g)\delta_{hg}$
for the real numbers $0\le c(n,h,g)\le1$ given below.
Thus
\begin{equation}
\label{eq-cnhgbound}
\|\Phi_n\circ\Psi_n(\lambda_h)-\lambda_h\|=1-\inf_{g\in G}c(n,h,g)
\end{equation}
(We write $g=g_1g_2\cdots g_{\ell(g)}$ as a reduced word in the $G_\iota$'s with the usual conventions explained above.)

\begin{case}
\label{case-PsiPhi-i}
$\ell(g)\le m$ and $\ell(hg)\le m$.
\end{case}
\noindent
Then $\Psi_n\circ\Phi_n(\lambda_h)\delta_g=V_n^*(\delta_{hg}\otimes\delta_e)=\delta_{hg}$ so $c(n,h,g)=1$.

\begin{case}
\label{case-PsiPhi-ii}
$\ell(g)\le m$ and $m\le\ell(hg)\le n$.
\end{case}
\noindent
Then $\Psi_n\circ\Phi_n(\lambda_h)\delta_g=V_n^*(\delta_{hg}\otimes\delta_e)=\frac{\sqrt{n-\ell(hg)}}{\sqrt{n-m}}\delta_{hg}$
and hence
\begin{equation}
\label{eq-cnhg-i}
\sqrt{1-\frac{\ell(h)}{n-m}}\le c(n,h,g)\le1.
\end{equation}

\begin{case}
\label{case-PsiPhi-iii}
$m\le\ell(g)\le n$ and $\ell(hg)\le m$.
\end{case}
\noindent
Then 
\begin{equation}
\label{eq-Vnstar-iii} 
\Psi_n\circ\Phi_n(\lambda_h)\delta_g=
V_n^*\left(\frac{\sqrt{n-\ell(g)}}{\sqrt{n-m}}\delta_{hg}\otimes\delta_e
+\sum_{j=m+\ell(g)-\ell(hg)}^{\ell(g)-1}\frac1{\sqrt{n-m}}\delta_{hg_1\ldots g_j}\otimes\delta_{g_{j+1}\ldots g_{\ell(g)}}\right),
\end{equation}
and the above quantity is equal to
$\frac{\sqrt{n-\ell(g)}}{\sqrt{n-m}}\delta_{hg}$;
hence~\eqref{eq-cnhg-i} holds.

\begin{case}
\label{case-PsiPhi-iv}
$m\le\ell(g)\le n$ and $m\le\ell(hg)\le n$.
\end{case}
\noindent
Then~\eqref{eq-Vnstar-iii} holds, whereas
\begin{equation}
\label{eq-Vn-iv}
V_n\delta_{hg}=
\frac{\sqrt{n-\ell(hg)}}{\sqrt{n-m}}\delta_{hg}\otimes\delta_e
+\sum_{j=m+\ell(g)-\ell(hg)}^{\ell(g)-1}\frac1{\sqrt{n-m}}\delta_{hg_1\ldots g_j}\otimes\delta_{g_{j+1}\ldots g_{\ell(g)}}.
\end{equation}
Hence
\begin{equation*}
\Psi_n\circ\Phi_n(\lambda_h)\delta_g=
\left(\frac{\sqrt{(n-\ell(g))(n-\ell(hg))}}{n-m}+\frac{\min(\ell(g),\ell(hg))-m}{n-m}\right)\delta_{hg}
\end{equation*}
and
\begin{equation*}
1-\frac{\ell(h)}{n-m}\le\frac{n-\max(\ell(g),\ell(hg))}{n-m}+\frac{\min(\ell(g),\ell(hg))-m}{n-m}
\le c(n,h,g)\le1.
\end{equation*}

\begin{case}
\label{case-PsiPhi-v}
$m\le\ell(g)\le n$ and $n\le\ell(hg)$.
\end{case}
\noindent
Then~\eqref{eq-Vnstar-iii} holds, whereas
\begin{equation}
\label{eq-Vn-v}
V_n\delta_{hg}=
\sum_{j=m+\ell(g)-\ell(hg)}^{n+\ell(g)-\ell(hg)-1}\frac1{\sqrt{n-m}}\delta_{hg_1\ldots g_j}\otimes\delta_{g_{j+1}\ldots g_{\ell(g)}}.
\end{equation}
Hence
\begin{equation*}
\Psi_n\circ\Phi_n(\lambda_h)\delta_g=
\frac{n+\ell(g)-\ell(hg)-m}{n-m}\delta_{hg}
\end{equation*}
and
\begin{equation}
\label{eq-cnhg-v}
1-\frac{\ell(h)}{n-m}\le c(n,h,g)\le1.
\end{equation}

\begin{case}
\label{case-PsiPhi-vi}
$n\le\ell(g)$ and $m\le\ell(hg)\le n$.
\end{case}
\noindent
Then
\begin{equation}
\label{eq-Vnstar-vi}
\Psi_n\circ\Phi_n(\lambda_h)\delta_g=
V_n^*\left(\sum_{j=m}^{n-1}\frac1{\sqrt{n-m}}\delta_{hg_1\ldots g_j}\otimes\delta_{g_{j+1}\ldots g_{\ell(g)}}\right)
\end{equation}
whereas~\eqref{eq-Vn-iv} holds.
Hence
\begin{equation*}
\Psi_n\circ\Phi_n(\lambda_h)\delta_g=
\frac{n-\ell(g)+\ell(hg)-m}{n-m}\delta_{hg}
\end{equation*}
and~\eqref{eq-cnhg-v} holds.

\begin{case}
\label{case-PsiPhi-vii}
$n\le\ell(g)$ and $n\le\ell(hg)$.
\end{case}
\noindent
Then~\eqref{eq-Vnstar-vi} and~\eqref{eq-Vn-v} hold.
Hence 
\begin{equation*}
\Psi_n\circ\Phi_n(\lambda_h)\delta_g=
\frac{n-m-|\ell(g)-\ell(hg)|}{n-m}\delta_{hg}
\end{equation*}
and~\eqref{eq-cnhg-v} holds.

Considering the bound~\eqref{eq-cnhgbound}
and the estimates in cases~\ref{case-PsiPhi-i}--\ref{case-PsiPhi-vii}, the proposition is proved.
\end{proof}

The following lemma will use Kirchberg's result~\cite[7.1]{Kirchberg:ComUHF} (see also~\cite[Prop.\ 2]{HRV:exact})
that if $0\to J\to B\to A\to 0$ is a short exact sequence
of C$^*$--algebras and $*$--homomorphisms, if this sequence has a completely positive contractive splitting
$s:A\to B$ and if $J$ and $A$ are exact C$^*$--algebras then $B$ is an exact C$^*$--algebra.

\begin{lemma}
\label{lem-stepII}
Suppose that every $G_\iota$ is infinite and $I$ is finite.
Then for every $n\in\Nats$ there is an exact C$^*$--algebra $D_n\subseteq\Bc(l^2(W_n))$ such that
$\Phi_n(C^*_r(G))\subseteq D_n$.
\end{lemma}
\begin{proof}
Since each $G_\iota$ is infinite, the reduced group C$^*$--algebra $C^*_r(G_\iota)$ contains
no nonzero compact operators on $l^2(G_\iota)$.
Write $G_\iota\oup=G_\iota\backslash\{e\}$ and let $Q_\iota\oup:l^2(G_\iota)\to l^2(G_\iota\oup)$
be the projection.
Let $\Ahat_\iota\subseteq\Bc(l^2(G_\iota\oup))$ be the C$^*$--algebra generated by
$\{Q_\iota\oup aQ_\iota\oup\mid a\in C^*_r(G_\iota)\}\cup\KEu(l^2(G_\iota\oup))$.
Since for every $a_1,a_2\in C^*_r(G_\iota)$ the operator
$Q_\iota\oup a_1a_2Q_\iota\oup-Q_\iota\oup a_1Q_\iota\oup a_2Q_\iota\oup$ has rank one,
we see that the quotient $\Ahat_\iota/\KEu(l^2(G_\iota))$ is isomorphic to $C^*_r(G_\iota)$.
Moreover, the short exact sequence 
\begin{equation}
\label{eq-Ahses}
0\to\KEu(l^2(G_\iota))\to\Ahat_\iota\to C^*_r(G_\iota)\to0
\end{equation}
has the unital completely positive splitting $a\mapsto Q_\iota\oup aQ_\iota\oup$.
By Kirchberg's result as described above, $\Ahat_\iota$ is an exact C$^*$--algebra.

Writing elements of $G$ as words in $(G_\iota\oup)_{\iota\in I}$, we canonically identify $l^2(W_n)$ with
\begin{equation}
\label{eq-l2Wntens}
\Cpx\delta_e\oplus\bigoplus_{\substack{1\le k\le n \\ \iota_1,\ldots,\iota_k\in I \\
  \iota_1\ne\iota_2,\ldots,\iota_{k-1}\ne\iota_k}}
l^2(G_{\iota_1}\oup)\otdt l^2(G_{\iota_k}\oup).
\end{equation}
Now we single out the tensor factors that are a fixed number of places counting from the right.
For $\iota\in I$ and $1\le p\le n$ consider the Hilbert spaces
\begin{align*}
\LEu_{p,\iota}&=\begin{cases}
  \displaystyle\Cpx\eta\oplus\bigoplus_{\substack{1\le k\le n-p-1 \\ \iota_1,\ldots,\iota_k\in I \\
        \iota_1\ne\iota_2,\ldots,\iota_{k-1}\ne\iota_k \\ \iota_k\ne\iota}}
    l^2(G_{\iota_1}\oup)\otdt l^2(G_{\iota_k}\oup)&\text{if }p<n \\[1ex]
  \hbox{\hspace{9em}}\Cpx\eta&\text{if }p=n \end{cases} \\[2ex]
\REu_{p,\iota}&=\begin{cases}
  \hbox{\hspace{9em}}\Cpx\eta&\text{if }p=1  \\[1ex]
  \displaystyle\bigoplus_{\substack{\iota_1,\ldots,\iota_{p-1}\in I \\
        \iota_1\ne\iota_2,\ldots,\iota_{p-2}\ne\iota_{p-1} \\ \iota_1\ne\iota}}
    l^2(G_{\iota_1}\oup)\otdt l^2(G_{\iota_{p-1}}\oup)&\text{if }p>1, \\[1ex] \end{cases}
\end{align*}
where $\Cpx\eta$ denotes a copy of $\Cpx$ having a unit vector $\eta$.
Let
\begin{equation*}
V_{p,\iota}:\LEu_{p,\iota}\otimes l^2(G_\iota\oup)\otimes\REu_{\iota,p}\to l^2(W_n)
\end{equation*}
be the isometric mapping of Hilbert spaces obtained by concatenating tensor products,
using the identification of $l^2(W_n)$ with~\eqref{eq-l2Wntens} and by
absorbing tensor factor copies of $\Cpx\eta$.
Then the range of $V_{p,\iota}$ is spanned by those words having an element of $G_\iota\oup$ in the $p$th position
counting from the right.
Let $J_0=\KEu(l^2(W_n))$ and for every $1\le p\le n$ consider the C$^*$--subalgebras of $\Bc(L^2(W_n))$,
\begin{equation}
\label{eq-Jpdef}
J_p=\bigoplus_{\iota\in I}V_{p,\iota}\bigl(\KEu(\LEu_{p,\iota})\otimes\Ahat_\iota\otimes1_{\REu_{p,\iota}}\bigr)V_{p,\iota}^*,
\end{equation}
where we mean that $J_p$ is the C$^*$--algebra generated by
$\bigcup_{\iota\in I}V_{p,\iota}\bigl(\KEu(\LEu_{p,\iota})\otimes\Ahat_\iota\otimes1_{\REu_{p,\iota}}\bigr)V_{p,\iota}^*$;
we have used the direct sum notation because the $V_{p,\iota}$ have orthogonal ranges.
Let $D_n=J_0+J_1+\cdots+J_n$.

We will now show that $D_n$ is an exact C$^*$--algebra.
The crucial observation is that $J_pJ_q\subseteq J_{\min(p,q)}$ whenever $0\le p,q\le n$.
Letting $B_p=J_0+J_1+\cdots+J_p$, we therefore have that $B_{p-1}$ is a two--sided ideal of $B_p$
for every $1\le p\le n$.
We will show by induction on $p\in\{0,1,\ldots,n\}$ that $B_p$ is an exact C$^*$--algebra.
If $p=0$ then $B_0=\KEu(l^2(W_n))$ is an exact C$^*$--algebra.
For the induction step, suppose $1\le p\le n$ and $B_{p-1}$ is an exact C$^*$--algebra,
and let us show that $B_p$ is an exact C$^*$--algebra.
Considering quotient $*$--algebras we have
\begin{equation}
\label{eq-quotients}
\frac{B_p}{B_{p-1}}=\frac{B_{p-1}+J_p}{B_{p-1}}\cong\frac{J_p}{B_{p-1}\cap J_p},
\end{equation}
and the canonical isomorphism above is isometric with respect to the quotient norms.
Therefore $B_p/B_{p-1}$ is a C$^*$--algebra,
and hence $B_p$ is a C$^*$--algebra.
We have the commutative diagram 
\begin{equation}
\label{eq-comd}
\begin{array}{ccccccccc}
0&\longrightarrow&B_{p-1}&\longrightarrow&B_p&\longrightarrow&B_p/B_{p-1}&\longrightarrow&0 \\[2ex]
&&\bigcup&&\bigcup&&\Big\uparrow\\[2ex]
0&\longrightarrow&B_{p-1}\cap J_p&\longrightarrow&J_p&\longrightarrow&J_p/(B_{p-1}\cap J_p)&\longrightarrow&0,
\end{array}
\end{equation}
where the two horizontal sequences are exact
and the upwards directed arrow is the isomporphism from~\eqref{eq-quotients}.
Recalling the definition of $J_p$ in~\eqref{eq-Jpdef}, we see that
\begin{equation*}
B_{p-1}\cap J_p=\bigoplus_{\iota\in I}
V_{p,\iota}\bigl(\KEu(\LEu_{p,\iota})\otimes\KEu(l^2(G_\iota\oup))\otimes1_{\REu_{p,\iota}}\bigr)V_{p,\iota}^*.
\end{equation*}
Hence
\begin{equation*}
\frac{J_p}{B_{p-1}\cap J_p}\cong\bigoplus_{\iota\in I}\KEu(\LEu_{p,\iota})\otimes\bigl(\Ahat_\iota/\KEu(l^2(G_\iota\oup))\bigr)
\cong\bigoplus_{\iota\in I}\KEu(\LEu_{p,\iota})\otimes C^*_r(G_\iota).
\end{equation*}
Using the hypothesis that each $G_\iota$ is exact, we have that $J_p/(B_{p-1}\cap J_p)$ is an exact C$^*$--algebra.
The above isomorphisms and the splittings found for the short exact sequences~\eqref{eq-Ahses} can be used to construct
a completely positive contractive splitting for the bottom horizontal sequence in~\eqref{eq-comd}.
This in turn gives a completely positive contractive splitting for the top horizontal sequence in~\eqref{eq-comd}.
Now appealing to Kirchberg's result as described before the statement of this lemma,
we conclude that $B_p$ is an exact C$^*$--algebra.
Therefore $D_n=B_n$ is an exact C$^*$--algebra.

It remains to show that $\Phi_n(C^*_r(G))\subseteq D_n$, and to do so we will use for the first time
the hypothesis that $I$ is finite.
Let us investigate $\Phi_n(\lambda_h)$ for $h\in G\backslash\{e\}$;
write $h=h_1h_2\cdots h_m$ where $h_j\in G_{\iota_j}\oup$, $\iota_j\ne\iota_{j+1}$.
Let $g\in W_n\backslash\{e\}$ and write $g=g_1g_2\cdots g_k$ for $g_j\in G_{\iota_j'}\oup$, $\iota_j'\ne\iota_{j+1}'$.
The value of $P_n\lambda_hP_n\delta_g$ depends on how much cancellation there is when multiplying $hg$.
On the right--hand--sides of the following equations, we will identify $l^2(W_n)$ with the Hilbert space~\eqref{eq-l2Wntens}.
\begin{case}
\label{case-Phinh-i}
$\iota_m\ne\iota_1'$.
\end{case}
\noindent
Then
\begin{equation*}
P_n\lambda_hP_n\delta_g=\begin{cases}
\delta_{h_1}\otdt\delta_{h_m}\otimes\delta_{g_1}\otdt\delta_{g_k}&\text{if }m+k\le n \\
0&\text{if }m+k>n.
\end{cases}
\end{equation*}

\begin{case}
\label{case-Phinh-ii}
$\iota_m=\iota_1',\,\iota_{m-1}=\iota_2',\ldots,\iota_{m-q+1}=\iota_q'$ but $\iota_{m-q}\ne\iota_{q+1}'$;
and $h_m=g_1^{-1},\,h_{m-1}=g_2^{-1},\ldots,h_{m-q+1}=g_q^{-1}$ for some $1\le q<\min(m,k)$.
\end{case}
\noindent
Then
\begin{equation*}
P_n\lambda_hP_n\delta_g=\begin{cases}
\delta_{h_1}\otdt\delta_{h_{m-q}}\otimes\delta_{g_{q+1}}\otdt\delta_{g_k}&\text{if }m+k-2q\le n \\
0&\text{if }m+k-2q>n.
\end{cases}
\end{equation*}

\begin{case}
\label{case-Phinh-iii}
$m<k$, $\iota_m=\iota_1',\ldots,\iota_1=\iota_m'$ and $h_m=g_1^{-1},\ldots,h_1=g_m^{-1}$.
\end{case}
\noindent
Then
\begin{equation*}
P_n\lambda_hP_n\delta_g=\delta_{g_{m+1}}\otdt\delta_{g_k}.
\end{equation*}

\begin{case}
\label{case-Phinh-iv}
$m>k$, $\iota_m=\iota_1',\ldots,\iota_{m-k+1}=\iota_k'$ and $h_m=g_1^{-1},\ldots,h_{m-k+1}=g_k^{-1}$.
\end{case}
\noindent
Then
\begin{equation*}
P_n\lambda_hP_n\delta_g=\begin{cases}
\delta_{h_1}\otdt\delta_{h_{m-k}}&\text{if }m-k\le n \\
0&\text{if }m-k>n.
\end{cases}
\end{equation*}
\pagebreak[1]
\begin{case}
\label{case-Phinh-v}
$m=k$, $\iota_m=\iota_1',\ldots,\iota_1=\iota_k'$ and $h_m=g_1^{-1},\ldots,h_1=g_k^{-1}$.
\end{case}
\noindent
Then
\begin{equation*}
P_n\lambda_hP_n\delta_g=\delta_e
\end{equation*}

\begin{case}
\label{case-Phinh-vi}
$\iota_m=\iota_1'$ and $h_m\ne g_1^{-1}$.
\end{case}
\noindent
Then
\begin{equation*}
P_n\lambda_hP_n\delta_g=\begin{cases}
\delta_{h_1}\otdt\delta_{h_{m-1}}\otimes\delta_{h_mg_1}\otimes\delta_{g_2}\otdt\delta_{g_k}&\text{if }m+k-1\le n \\
0&\text{if }m+k-1>n.
\end{cases}
\end{equation*}

\begin{case}
\label{case-Phinh-vii}
$\iota_m=\iota_1',\,\iota_{m-1}=\iota_2',\ldots,\iota_{m-q+1}=\iota_q',\,\iota_{m-q}'=\iota_{q+1}'$;
and $h_m=g_1^{-1},\,h_{m-1}=g_2^{-1},\ldots,h_{m-q+1}=g_q^{-1}$ but $h_{m-q}\ne g_{q+1}^{-1}$
for some $1\le q<\min(m,k)$.
\end{case}
\noindent
Then
\begin{equation*}
P_n\lambda_hP_n\delta_g=\begin{cases}
\delta_{h_1}\otdt\delta_{h_{m-q-1}}\otimes\delta_{h_{m-q}g_{q+1}}\otimes\delta_{g_{q+2}}\otdt\delta_{g_k}&\text{if }m+k-2q-1\le n \\
0&\text{if }m+k-2q-1>n,
\end{cases}
\end{equation*}
where in the above equation if $q=k-1$ then the expression $\otimes\delta_{g_{q+2}}\otdt\delta_{g_k}$ should be treated as absent,
and if $q=m-1$ then the expression $\delta_{h_1}\otdt\delta_{h_{m-q-1}}\otimes$ should be treated as absent.

Cases~\ref{case-Phinh-i}--\ref{case-Phinh-vii} are the only possibilities;
using them we get the following expression:
\begin{align}
\Phi_n(\lambda_h)=&
\displaystyle\sum_{r=\max(m-n,0)}^{\min(m,n)}
 \langle\delta_{h_1}\otdt\delta_{h_r},\delta_{h_m^{-1}}\otdt\delta_{h_{r+1}^{-1}}\rangle_\KEu \label{eq-Phinlh-i} \\[2ex]
&\begin{aligned}
  \displaystyle+\sum_{r=\max(m-n+1,0)}^{\min(m,n-1)}&\sum_{p=1}^{n-\max(r,m-r)}\sum_{\iota\in I\backslash\{\iota_r,\iota_{r+1}\}}
   V_{p,\iota}\biggl( \\
  &\bigl\langle\delta_{h_1}\otdt\delta_{h_r},\delta_{h_m^{-1}}\otdt\delta_{h_{r+1}^{-1}}\bigr\rangle_\KEu
   \otimes1_{l^2(G_\iota\oup)}\otimes1_{\REu_{p,\iota}}\biggr)V_{p,\iota}^* \label{eq-Phinlh-ii}\end{aligned} \\[2ex]
&\begin{aligned}
  \displaystyle+\sum_{r=\max(m-n,0)+1}^{\min(m,n)}&\sum_{p=1}^{n-\max(r-1,m-r)}V_{p,\iota_r}\biggl( \\
  &\bigl\langle\delta_{h_1}\otdt\delta_{h_{r-1}},\delta_{h_m^{-1}}\otdt\delta_{h_{r+1}^{-1}}\bigr\rangle_\KEu\otimes
   Q_{\iota_r}\oup\lambda_{h_r}Q_{\iota_r}\oup\otimes1_{\REu_{p,\iota_r}}\biggr)V_{p,\iota_r}^*,\end{aligned} \label{eq-Phinlh-iii}
\end{align}
where we use the conventions that if $v$ and $w$ are vectors in a Hilbert space then
$\langle v,w\rangle_\KEu$ denotes the rank one operator given by
$\bigl(\langle v,w\rangle_\KEu\bigr)(\zeta)=v\langle w,\zeta\rangle$, with the inner product on the right being in the Hilbert space and
linear in the second variable, and where the symbols
\begin{alignat*}{2}
\delta_{h_1}\otdt\delta_{h_r}&\text{ in line \eqref{eq-Phinlh-i} }&\text{should be interpreted to mean }\delta_e&\text{ if }r=0
\displaybreak[0] \\
\delta_{h_m^{-1}}\otdt\delta_{h_{r+1}^{-1}}&\text{ in line \eqref{eq-Phinlh-i} }&\text{should be interpreted to mean }\delta_e&\text{ if }r=m
\displaybreak[0] \\
\delta_{h_1}\otdt\delta_{h_r}&\text{ in line \eqref{eq-Phinlh-ii} }&\text{should be interpreted to mean }\eta&\text{ if }r=0
\displaybreak[0] \\
\delta_{h_m^{-1}}\otdt\delta_{h_{r+1}^{-1}}&\text{ in line \eqref{eq-Phinlh-ii} }&\text{should be interpreted to mean }\eta&\text{ if }r=m
\displaybreak[0] \\
\delta_{h_1}\otdt\delta_{h_{r-1}}&\text{ in line \eqref{eq-Phinlh-iii} }&\text{should be interpreted to mean }\eta&\text{ if }r=1
\displaybreak[0] \\
\delta_{h_m^{-1}}\otdt\delta_{h_{r+1}^{-1}}&\text{ in line \eqref{eq-Phinlh-iii} }&\text{should be interpreted to mean }\eta&\text{ if }r=m.
\end{alignat*}
Note that in the above formula for $\Phi_n(\lambda_h)$,
cases~\ref{case-Phinh-iii} and~\ref{case-Phinh-iv} as well as the case $g=e$ are incorporated into line~\eqref{eq-Phinlh-i},
cases~\ref{case-Phinh-i}, \ref{case-Phinh-ii} and~\ref{case-Phinh-iii} are incorporated into line~\eqref{eq-Phinlh-ii}
and cases~\ref{case-Phinh-vi} and~\ref{case-Phinh-vii} are incorporated into line~\eqref{eq-Phinlh-iii}.

Thus $\Phi_n(\lambda_h)\in D_n$;
hence $\Phi_n(C^*_r(G))\subseteq D_n$.
\end{proof}

\begin{thm}
Let $I$ be a nonempty set and for every $\iota\in I$ let $G_\iota$ be an exact group.
Let $G=\freeprod_{\iota\in I} G_\iota$ be the free product of groups.
Then $G$ is an exact group.
\end{thm}
\begin{proof}
Since $G$ is the direct limit of free products of finite subfamilies of $(G_\iota)_{\iota\in I}$, we may
without loss of generality assume that $I$ is finite.
Furthermore, we may without loss of generality assume
that each $G_\iota$ is infinite, by replacing $G_\iota$ with the exact group $G_\iota\oplus\Ints$ if necessary
and then passing to the subgroup of the free product.
Using Lemmas~\ref{lem-stepI} and~\ref{lem-stepII} we find for every $n\ge1$ an exact C$^*$--algebra $D_n$ and
completely positive contractive maps $\Phi_n:C^*_r(G)\to D_n$ and $\Psi_n:D_n\to\Bc(l^2(G))$ such that
$\lim_{n\to\infty}\|\Psi_n\circ\Phi_n(a)-a\|=0$ for every $a\in C^*_r(G)$.
Now Lemma~\ref{lem-W} shows that $C^*_r(G)$ is exact.
\end{proof}

\end{spacing}
\begin{spacing}{1.0}

\bibliographystyle{plain}

\begin{thebibliography}{99}

\bibitem{Adams:amen} S.\ Adams,
{\em Boundary amenability for word hyperbolic groups and an application to smooth dynamics of simple groups},
Topology {\bf 33} (1994), 765-783.

\bibitem{Choi:Z2Z3} M.-D.\ Choi,
{\em A simple C$^*$--algebra genereated by two finite--order unitaries},
Canad.\ J.\ Math.\ {\bf 31} (1979), 867-880.

\bibitem{D:ExactFP} K.J.~Dykema,
{\em Exactness of reduced amalgamated free products of {$C^*$}-algebras,}
preprint (1999).

\bibitem{Ger:exgp} E.\ Germain,
{\em Approximate invariant means for boundary actions of hyperbolic groups},
appendix to {\em Amenable Groupoids} by C.\ Anatharaman--Delaroche and J.\ Renault,
preprint.

\bibitem{HRV:exact} P.\ de la Harpe, A.G.\ Robertson, A.\ Valette,
{\em On exactness of group C$^*$--algebras,}
Quart.\ J.\ Math.\ Oxford (2) {\bf 45} (1994), 499-513.

\bibitem{Kirchberg:ComUHF} E.\ Kirchberg,
{\em Commutants of unitaries in UHF algebras and functorial properties of exactness,}
J.\ reine angew.\ Math.\ {\bf 452} (1994), 39-77.

\bibitem{KW:exgp} E.\ Kirchberg, S.\ Wasermann,
{\em Exact groups and continuous bundles of C$^*$--algebras},
preprint.

\bibitem{KW:perm} \rule{3em}{.1mm},
{\em Permanence properties of C$^*$--exact groups},
preprint.

\bibitem{Wa76} S.\ Wassermann,
{\em On tensor products of certain group C$^*$--algebras},
J.\ Funct.\ Anal.\ {\bf 23} (1976), 239-254.

\bibitem{Wa:nuclemb} \rule{3em}{.1mm},
{\em Tensor products of free--group C$^*$--algebras},
Bull.\ London Math.\ Soc.\ {\bf 22} (1990), 375-380.

\bibitem{Wa:SNU} \rule{3em}{.1mm},
{\em Exact C$^*$--algebras and Related Topics},
Seoul National University Lecture Notes Series {\bf 19}, 1994.

\end{thebibliography}

\vskip1ex
\scriptsize
\noindent{\sc Department of Mathamtics, Texas A\&M University, College Station TX 77843--3368, USA} \\
\noindent{\sl E-mail:} {\tt Ken.Dykema@math.tamu.edu} \\
\noindent{\sl Internet URL:} {\tt http://www.math.tamu.edu/$\,\widetilde{\;}$Ken.Dykema/} \\

\end{spacing}
\end{document}